\documentclass{amsart}
\usepackage{amsmath}
\usepackage{amssymb}
\usepackage{amsthm}
\usepackage{url}
\newtheorem{theorem}{Theorem}[section]

\newtheorem*{axiom}{External induction}

\newtheorem{definition}[theorem]{Definition}

\newcommand*{\st}{\mathrm{st}}
\newcommand*{\forallst}{\forall^\mathrm{st}}

\begin{document}

\title{Old and new approaches to the Sorites paradox}
\author{Bruno Dinis}

\thanks{The author would like to thank Imme van den Berg and Pedro Isidoro for several remarks on preliminary versions of this paper. }


\begin{abstract}
  The Sorites paradox is the name of a class of paradoxes that arise when vague predicates are considered. Vague predicates lack sharp boundaries in extension and is therefore not clear exactly when such predicates apply.  
   Several approaches to this class of paradoxes have been made since its first formulation by Eubulides of Miletus in the IV century BCE. In this paper I survey some of these approaches and point out some of the criticism that these approaches have received. A new approach that uses tools from nonstandard analysis to model the paradox is proposed. 

\end{abstract}

\maketitle

\begin{quotation}
The knowledge came upon me, not quickly, but little by little and grain by
grain.

(Charles Dickens in David Copperfield)
\end{quotation}

\section{To be or not to be a heap}
This paper concerns the paradoxes which arise when several orders of mag-
nitude are considered. These paradoxes are part of the larger phenomenon known as
vagueness. Indeed,

\begin{quote}
[...] vagueness can be seen to be a feature of syntactic categories other than predicates. Names, adjectives, adverbs and so on are all susceptible to paradoxical sorites reasoning in a derivative sense. \cite{Hyde Dominic}
\end{quote}

Vague predicates share at least the following features:

\begin{enumerate}
\item Admit borderline cases

\item Lack sharp boundaries

\item Are susceptible to sorites paradoxes
\end{enumerate}

\emph{Borderline} cases are the ones where it is not clear whether or not
the predicate applies, independently of how much one knows about it. For
instance, most basketball players are clearly tall and most jockeys are
clearly not tall. But in many cases is rather unclear if the person in
question is tall, even if one knows its height with great precision.
Furthermore, there is no clear distinction between the set of all tall
people and the set of people that are not tall. These sets lack sharp
boundaries. This leads to a collection of paradoxes called \emph{Sorites}
paradoxes first formulated by the Ancient Greek philosopher Eubulides of Miletus in the IV century BCE. One can be stated in the following way: a single grain of wheat
cannot be considered as a heap. Neither can two grains of wheat. One must
admit the presence of a heap sooner or later, so where to draw the line? In
fact, the name Sorites derives from the Greek word \emph{sor\'os} which means heap. However, one
can reconstruct the paradox by replacing the term 'heap' by other vague
concepts such as 'tall', 'beautiful', 'bald', 'heavy', 'cold', 'rich',...

The argument consists of a predicate $S$ (the soritical predicate) and a
subject expression $a_{n}$\ in the series regarding to which $S$ is
soritical. The terms of the series are supposed to be ordered. According to
Barnes \cite{Barnes} a predicate $S$ must satisfy three constraints in order
to be considered soritical:

\begin{enumerate}
\item Appear to be valid for $a_{1}$, the first item in the series;

\item Appear to be false for $a_{i}$, the last item in the series;

\item Each adjacent pair in the series, $a_{n}$ and $a_{n+1}$ must be
sufficiently similar as to appear indiscriminable in respect to $S$.
\end{enumerate}

This means that the predicate $S$ needs to be sufficiently vague in order to
allow small changes. Small changes do not determine the difference between a
set of individual grains and a heap, between a bald man and a hairy one,
between a rich person and a poor one. However, and in spite of the vagueness
involved, it also needs to have a certain area on which $S$ is clearly true
and an area on which $S$ is clearly false.

\begin{quotation}
The difference of one grain would seem to be too small to make any
difference to the application of the predicate; it is a difference so
negligible as to make no apparent difference to the truth-values of the
respective antecedents and consequents.

Yet the conclusion seems false. \cite{Hyde Dominic}
\end{quotation}

 This paper surveys some approaches that have been made to deal with the phenomenon of vagueness and the Sorites type paradoxes which arise when vague predicates are used. Also, a new approach that uses tools from nonstandard analysis to model the paradox is proposed. As mentioned above, soritical arguments are tolerant to small changes but not tolerant to large changes in relevant aspects. In fact, using a special class of external sets called external numbers (see Section \ref{EN}) it is possible to define rigorously what is meant with terms such as 'small changes' or 'large changes'. The fact that large changes come as the result of the accumulation of small changes is no surprise because it is a very well known fact from nonstandard analysis that an infinitely large sum of infinitesimals may very well become appreciable or even infinitely large. In this way, one can make a rigorous claim that a heap and a set of individual grains of wheat are indeed not of the same order of magnitude.

\section{Forms of the paradox}\label{section Paradoxical forms}

The Sorites paradox can be stated in various ways. This implies that one
cannot hope to solve the paradox by pointing out a fault particular to any
one of those way. One should instead try and reveal a common fault to all
possible forms that the paradox can take. I do so by considering the (standard)
mathematical induction and conditional schemata. These schemata will be revisited, in Section \ref{section nsa point of view}, after a nonstandard point of view is adopted. 
\subsection{Induction\label{section induction}}

Mathematical induction is generally used (within standard mathematics) to
prove that a mathematical statement involving a natural number $n$ holds for
all possible values of $n$. This is done in two steps. On the first step (%
\emph{basis}) one proves that there is a first element for which the
statement holds. On the second step (\emph{inductive step}) one shows that
if the statement holds for some $n$ then it also holds for $n+1$. Then, by the principle of mathematical induction, the
statement is valid for all $n$.
Let $S$ represent a soritical predicate, for example 'is not a
heap' and let $a_{n}$ represent the $n$-th element in the soritical series. In the example above it would be the sentence '$n$ grains of wheat'.

The Sorites paradox can now be represented in the following way:

\begin{equation*}
\left \{ 
\begin{array}{c}
\left( Sa_{1}\wedge \forall n\left( Sa_{n}\rightarrow Sa_{n+1}\right)
\right) \rightarrow \forall nSa_{n} \\ 
\exists \omega \left( \lnot Sa_{\omega }\right)%
\end{array}%
\right.
\end{equation*}%
So, if one admits that:

\begin{enumerate}
\item A single grain of wheat is not a heap.

\item If a collection of $n$ grains of wheat is not a heap then a collection
of $n+1$ grains of wheat is also not a heap.
\end{enumerate}

One concludes (by induction) that the heap will never appear. Since at some
point the heap is obviously there one might come to the conclusion that
there is something wrong with induction or, at least, with applying
induction to vague predicates.

\subsection{Conditional form\label{section conditional form}}

The conditional form of the Sorites paradox is the most common form
throughout the literature. Using the notation of the previous section it can be formalized in the following way:

If%
\begin{equation*}
\left \{ 
\begin{array}{c}
Sa_{1} \\ 
Sa_{1}\rightarrow Sa_{2} \\ 
Sa_{2}\rightarrow Sa_{3} \\ 
... \\ 
Sa_{i}\rightarrow Sa_{i+1} \\ 
\exists j\left( \lnot Sa_{j}\right)%
\end{array}%
\right.
\end{equation*}%
Assuming $Sa_{1}$, $Sa_{1}\rightarrow Sa_{2}$, $Sa_{2}\rightarrow Sa_{3},$
..., $Sa_{i}\rightarrow Sa_{i+1}$, by \textit{modus ponens, }the conclusion
is $Sa_{i}$, where $i$ can be arbitrarily large. This is a fairly simple
reasoning where the premises are: a single grain of wheat does not make a
heap; if one grain of wheat does not make a heap then two grains of wheat do
not form a heap either; if two grains of wheat do not make a heap then three
grains of wheat do not form a heap either... if $i$ grains of wheat do not
make a heap then $i+1$ grains of wheat do not form a heap either. The
conclusion is that a set of an arbitrarily large number of grains $i$ does
not make a heap. However if one observes that there is a set of $j$ grains
that form a heap it generates a paradox.\footnote{A reasoning from classical logic was used here and one could think that the paradox may be circumveined if one would use intuitionistic logic instead. In fact, there have been attempts to use intuitionistic logic to deal with the paradox (see \cite{Putnam} and its defense \cite{Schwartz}). However, according to Keefe, \begin{quote} [...] critics have shown that with various reasonable additional assumptions, other
versions of sorites arguments still lead to paradox. [...] The bulk of the criticisms point to the conclusion that
there is no sustainable account of vagueness that emerges from
rejecting classical logic in favour of intuitionistic logic.\cite[p. 22]{Keefe} \end{quote} For the criticism that Keefe refers to the reader may consult \cite{ReadWright,Chambers}.}

\section{Response attempts}

There are several attempts to solve the Sorites paradox. These responses are
divided into the following four types (\cite{Keefe}, p. 19-20). A first
type of response would be to deny the validity of the argument, refusing to
grant that the conclusion follows from the premises. Alternatively one can
question the strict truth of the inductive premise (or of one of the
conditionals). A third possibility is to accept the validity of the argument
and the truth of its inductive premise (or of all the conditional premises)
but contest the truth of the conclusion. Finally one can grant that there
are good reasons to consider both the argument form as valid, and accept the
premises and deny the conclusion hence proving that the predicate is
incoherent.

In this section, some of the responses to the paradox will be reviewed. For a wider
account on this matter see for example \cite{Keefe,Smith,Shapiro,Williamson}. I would like to emphasize that the theories presented
below correspond to a wide variety of related points of view. This means
that there are many versions of the theories presented. So, when reviewing a
theory I tend to give only the general lines, common to the various
versions of that theory.

\subsection{Ideal Languages}

Natural languages such as English or Vietnamese distinguish between
intension and extension of terms. The \emph{intension} is the internal
content of a term or concept while the \emph{extension} is the range of
applicability of a term by naming the particular objects that the term
denotes. The two predicates 'is a creature with a heart' and 'is a creature
with a kidney' (see \cite{Quine}) have the same extension because the set of
creatures with hearts and the set of creatures with kidneys are the same.
However, having a heart and having a kidney are very different things, so
one concludes that

\begin{quote}
terms can name the same thing but differ in meaning. \cite{Quine}
\end{quote}

Hence, the distinction between intension and extension leads necessarily to
vagueness, ambiguity and indeterminacy of meaning for words and phrases.
This is in part the reason why natural languages are so powerful. In poetry, for example, beauty is achieved by taking advantage of these features. However, if one needs clarity and precision of language then she is forced to conclude that natural
languages are not the way to go. According to Quine

\begin{quote}
The sorites paradox is one imperative reason for precision in science, along
with more familiar reasons. \cite{Quine1981}
\end{quote}

An ideal language would left out all such factors in order to eliminate any
vagueness.

The defenders of this response, among them Frege \cite{Frege}, Russell \cite%
{Russell} and Wittgenstein \cite{Wittgenstein PI}, consider vagueness as a
non-eliminable feature of natural language. The way to avoid vagueness is by creating and using ideal languages instead. This would mean that arguments of the Sorites type are not valid since they contain vague expressions. 


As stated by Russell,

\begin{quotation}
The fact is that all words are attributable without doubt over a certain
area, but become questionable within a penumbra, outside which they are
again certainly not attributable. Someone might seek to obtain precision in
the use of words by saying that no word is to be applied in the penumbra,
but unfortunately the penumbra is itself not accurately definable, and all
the vagueness which apply the primary use of words apply also when we try to
fix a limit to their indubitable applicability. \cite{Russell}
\end{quotation}

So, this response implies that it is the philosopher's job to discover a
logically ideal language. However, this doesn't seem possible using
classical logic:

\begin{quotation}
All traditional logic habitually assumes that precise symbols are being
employed. It is therefore not applicable to this terrestrial life, but only
to an imagined celestial existence. \cite{Russell}
\end{quotation}

Russell also believed that

\begin{quotation}
Vagueness, clearly, is a matter of degree, depending upon the extent of the
possible differences between different systems represented by the same
representation. Accuracy, on the contrary, is an ideal limit. \cite{Russell}
\end{quotation}

Ideal languages as a response to the sorites
paradox seem to have unsatisfying features. According to Keefe,

\begin{quotation}
denying the validity of the sorites argument seems to require giving up
absolutely fundamental rules of inference. \cite[p. 20]{Keefe}
\end{quotation}

So, if one chooses to go in this direction fundamental rules such as \textit{%
modus ponens }or mathematical induction are to be put in question. Furthermore, by eliminating vague predicates from the language one is not solving the paradox but avoiding it by sweeping it under the rug. In fact,
most philosophers nowadays believe that vagueness is an important part of
natural language and cannot be separated from it.

\subsection{The Epistemic theory} \label{Epistemic}

The Epistemic theory is based on the idea that the precise boundaries to
knowledge itself cannot be known. Vagueness is seen as a particular type of
ignorance.

The fact that this theory is built in the classical logic framework implies
that there are precise bounds for the extensions of vague predicates even if
one does not know where they are located. For instance, the defenders of the
epistemic theory claim that there is in fact a last grain of wheat in the
series before the heap turns up, even if one is not (nor ever will be) able
to identify it definitively. In fact, Williamson \cite{Williamson} has shown
that if there is a precise boundary for penumbral cases one cannot know
where it is. So, soritical predicates are indeterminate in extension but not
semantically. This position has been notably defended
by Williamson \cite{Williamson,Williamson2000} and Sorensen \cite%
{Sorensen,Sorensen2001}.

The first and major objection to this theory is its counter-intuitive
nature. The meaning of a word is (usually) determined by its use. According
to Wittgenstein

\begin{quote}
For a large class of cases - though not for all - in which we employ the
word 'meaning' it can be defined thus: the meaning of a word is its use in
the language \cite{Wittgenstein PI}
\end{quote}

and

\begin{quote}
if we had to name anything which is the life of the sign, we should have to
say that it is its use. \cite{Wittgenstein BB}
\end{quote}

For instance, the word 'guitar' means an actual guitar because one use that
word to mean an actual guitar (even if one does not know how to play). Now,
one does not usually use the word 'heap' as if a single grain of wheat could
make a difference. Neither, more generally, does one use any vague term as
if it were not tolerant to small changes. One does not use vague terms as if
they had precise borders. In this sense, Smith \cite{Smith} claims that the
epistemicist is forced to deny a link between meaning and use.

Another point that deserves criticism is that nothing is said about how
predicates get the precise extensions that they do. It is claimed that there
is in fact a last grain of wheat in the series before the heap turns up. So
there should be attempts to find which one is it \cite{Keefe}. I agree that
ignorance is no excuse for the lack of attempts to find the precise
boundaries of vague concepts. There should be at least some reasons to
believe about where these boundaries are.

\subsection{Supervaluationism}

According to Fine, vagueness is a semantic notion not to be confused with
ambiguity nor undecidability:

\begin{quotation}
Let us say, in a preliminary way, what vagueness is. I take it to be a
semantic notion. Very roughly, vagueness is deficiency of meaning. As such,
it is to be distinguished from generality, undecidability, and ambiguity.
These latter are, if you like, lack of content, possible knowledge, and
univocal meaning, respectively. \cite{Fine}
\end{quotation}

Supervaluationism proposes to solve the problem of vagueness by modifying
classical semantics, using Van Fraassen's supervaluations. According to Van
Fraassen:

\begin{quote}
A \textit{supervaluation over a model} is a function that assigns T (F)
exactly to those statements assigned T (F) by all the classical valuations
over that model. \cite{van Fraassen}
\end{quote}

And he concludes that

\begin{quote}
Supervaluations have truth-value gaps. \cite{van Fraassen}
\end{quote}

In classical logic the connectives have truth values in a functional way. I
recall that a connective of statements is \emph{truth-functional} if and
only if the truth value of any compound statement obtained by applying that
connective is a function of the individual truth values of the constituent
statements that form the compound. The classical logic connectives are all
truth-functional.\footnote{This is immediately visible if one computes the logical value of a given sentence using truth tables or a proof calculus like natural deduction \cite{Prawitz}.} Supervaluationists abandon the concept of
truth-functionality.

Fine applies the distinction between extension and intension \cite{Quine} to
vagueness:

\begin{quote}
Extensional vagueness is deficiency of extension, intensional vagueness
deficiency of intension. Moreover, if intension is the possibility of
extension, then intensional vagueness is the possibility of extensional
vagueness. \cite{Fine}
\end{quote}

According to this theory, a vague predicate does not need to have a unique,
sharply bounded, truth function. Vague predicates have things to which they
definitely apply (positive extension), things to which they definitely do
not (negative extension) and a penumbra (penumbral connections). The
penumbra involves cases which seem to be neither true nor false\footnote{%
Fine warns about the general confusion of under- and over-determinacy.
\par
\begin{quotation}
A vague sentence can be made more precise; and this operation should
preserve truth-value. But a vague sentence can be made to be either true or
false, and therefore the original sentence can be neither. \cite{Fine}
\end{quotation}
} (borderline cases). These penumbral connections are instances of
truth-value gaps. Truth-value gaps are related with extensional vagueness.
However,

\begin{quote}
Despite the connection, extensional vagueness should not be defined in terms
of truth-value gaps. This is because gaps can have other sources, such as
failure of reference or presupposition. \cite{Fine}
\end{quote}

Supervaluationists claim, roughly speaking, that a vague sentence is true if
and only if it is true for all ways of making it completely precise \cite%
{Fine}, called \emph{precisifications}. There are then many interpretations
or precisifications. Each one of these precisifications has no penumbra
because it behaves according to classical bivalence. The assignment of truth
value for all such precisifications is a supervaluation.

A sentence which is true in all precisifications is called \textit{\emph{%
supertrue }}and a sentence which is false in all precisifications is called 
\textit{\emph{superfalse}}. A sentence which is true for some
precisifications and false on others is neither true nor false\footnote{%
In fact, not all supervaluationists accept this last sentence.}. This means
in particular that tautologies from classical logic are supertrue.

According to Keefe \cite{Keefe}, truth is supertruth,  meaning that a sentence is true if and only if it is true on all admissible
precisifications. A precisification is acceptable only if the extensions of
the concepts do not overlap. The truth of a compound sentence is determined
by its truth on every precisification. For a wider account on supervaluationism the reader is referred to \cite{van Fraassen,Fine,Keefe}.

Fodor and Lepore are particularly critic of the supervaluationist approach
to vagueness:

\begin{quotation}
[...] there is something fundamentally wrong with using supervaluation
techniques either for preserving classical logic or for providing a
semantics for linguistic expressions ordinarily thought to produce
truth-value gaps. \cite{Fodor Lepore}
\end{quotation}

However the fault they point out is not of a logical nature. Indeed they say:

\begin{quotation}
Right from the start, however, we want to emphasize that the objections we
are raising are \textit{philosophical} rather than \textit{logical}. We have
no argument with supervaluations considered as a piece of formal
mathematics. \cite{Fodor Lepore}
\end{quotation}

Fodor and Lepore point out as the main flaws of supervaluationism the
violation of intuitive semantic principles concerning disjunctions and
existential quantification, the abandonment of classical rules of inference
and the violation of core principles concerning the concept of truth.

Let $S$ represent the predicate 'is not a heap'. Since all tautologies are
supertrue,%
\begin{equation*}
\lnot \left( \forall n\left( Sa_{n}\rightarrow Sa_{n+1}\right) \right)
\end{equation*}%
is equivalent to%
\begin{equation*}
\exists n\left( Sa_{n}\wedge \lnot Sa_{n+1}\right) \text{,}
\end{equation*}%
which, semantically speaking, seems to postulate the existence of a sharp
boundary and looks for that matter like a step back towards the epistemic
theory. Also $\left( S\vee \lnot S\right) $ is supertrue. So, for all
precisifications one of the statements is true. However, the statement $S$
is borderline and therefore neither true nor false.

Keefe \cite{Keefe} argues that it is possible to surpass these difficulties
at the price of adding a new operator to the language: the 'definitely'
operator $D$. This operator is however not closed under certain operations
such as contraposition and conditional introduction. So alternatives to the
classical closure principles are proposed. However, this implies that the
logic used is no longer classical.

Another argument against supervaluationism is that little information is
given on what makes a precisification acceptable other than saying that
precisifications must respect penumbral connections and therefore
their admissibility is a vague matter. Also, supervaluationism states that
precisifications behave in a classical way and have no penumbra. However,
each precisification may divide the positive and negative extensions in
different places. 

For more on objections to supervaluationism and attempts to respond to those objections the reader is referred to \cite{Cobreros}.

\subsection{Many-valued logics}

Many-valued logics is a general term that refers to logics which have more
than two truth-values. In these logics the principle of truth-functionality
is accepted and so a sentence remains unaffected when one of its components
is replaced by another with the same truth value. Many-valued logics became
accepted as an independent part of logic with the works of\  \L ukasiewicz
and Post in the 1920's. Since then many many-valued logics emerged (e.g. 
\cite{Lukasiewicz,Godel,Hajek et al.,Post}) and it is
not possible nor desirable to describe them all in these pages. However I
shall discuss an application of Kleene's three-valued logic and applications
of fuzzy logics because these seem to be the most relevant in what concerns the
phenomenon of vagueness. For a more complete reference concerning
many-valued logics see for example \cite{Gottwald}.

\subsubsection{Kleene's three-valued logic}

Perhaps one of the simplest and best-known examples of a many-valued logic
is Kleene's three-valued logic \cite{Kleene}. Kleene thought of the third
truth value as undefined or underdetermined\footnote{%
Priest \cite{Priest} gave an alternative three-valued logic conceiving the
third truth-value as overdetermined, interpreting the symbol $\frac{1}{2}$
as being both true and false.}. So one has three truth-values: $1$ (true), $%
0 $ (false) and $\frac{1}{2}$ (undefined or unknown). One has truth-tables
for which the connectives are \emph{regular}, i.e. in terms of ordering,
undefined is placed below both true and false. This means that the behavior
of the third truth value should be compatible with any increase in
information. Kleene proposed the following truth-tables for the so-called 
\emph{strong} connectives:

\begin{center}
\bigskip 
\begin{tabular}{|c|c|}
\hline
$p$ & $\lnot p$ \\ \hline
$1$ & $0$ \\ \hline
$0$ & $1$ \\ \hline
$\frac{1}{2}$ & $\frac{1}{2}$ \\ \hline
\end{tabular}
\  \  \ 
\begin{tabular}{|c|c|c|c|c|c|}
\hline
$p$ & $q$ & $p\vee q$ & $p\wedge q$ & $p\rightarrow q$ & $%
p\longleftrightarrow q$ \\ \hline
$1$ & $1$ & $1$ & $1$ & $1$ & $1$ \\ \hline
$1$ & $0$ & $1$ & $0$ & $0$ & $0$ \\ \hline
$1$ & $\frac{1}{2}$ & $1$ & $\frac{1}{2}$ & $\frac{1}{2}$ & $\frac{1}{2}$ \\ 
\hline
$0$ & $1$ & $1$ & $0$ & $1$ & $0$ \\ \hline
$0$ & $0$ & $0$ & $0$ & $1$ & $1$ \\ \hline
$0$ & $\frac{1}{2}$ & $\frac{1}{2}$ & $0$ & $1$ & $\frac{1}{2}$ \\ \hline
$\frac{1}{2}$ & $1$ & $1$ & $\frac{1}{2}$ & $1$ & $\frac{1}{2}$ \\ \hline
$\frac{1}{2}$ & $0$ & $\frac{1}{2}$ & $0$ & $\frac{1}{2}$ & $\frac{1}{2}$ \\ 
\hline
$\frac{1}{2}$ & $\frac{1}{2}$ & $\frac{1}{2}$ & $\frac{1}{2}$ & $\frac{1}{2}$
& $\frac{1}{2}$ \\ \hline
\end{tabular}%
\bigskip
\end{center}

These tables are uniquely determined as the strongest possible regular
extensions of the classical two-valued tables. Quantifiers can be defined in
the following way: $\exists x:P\left( x\right) $ is true if $P\left(
x\right) $ is true for some value of $x$ and it is false if $P\left(
x\right) $ is false for all values and indefinite otherwise; $\forall
xP\left( x\right) $ is true if $P\left( x\right) $ is true for all values of 
$x$ and false if $P\left( x\right) $ is false for some value and indefinite
otherwise. Tye \cite{Tye} applies Kleene's three-valued logic to the sorites
paradox. However, the objections made to the bipartite division can also be
used to refute a tripartite division. In fact, Tye \cite{Tye} claims that

\begin{quote}
[...] vagueness cannot be reconciled with any precise dividing lines.
\end{quote}

because

\begin{quote}
there is no determinate fact of the matter about where truth-value changes
occur.
\end{quote}

That is to say that there is no way to assign precise truth-values to vague
terms. So, as a solution, Tye proposes to use a vague metalanguage. He
claims that there are sets which are genuinely vague items. For instance the
set of tall men has borderline members (men which are neither clearly
members nor clearly non-members of the set).

\begin{quote}
There is no determinate fact of the matter about there are objects that are
neither members, borderline members, nor non-members. \cite{Tye}
\end{quote}

Kleene's three-valued logic has the undesirable feature of having no
tautologies, because the two-valued tautologies can take the value $\frac{1}{%
2}$ in the three-valued case. As an example consider the law of excluded
middle $p\vee \lnot p$. In Kleene's three-valued logic the truth table is
the following:

\begin{center}
\bigskip 
\begin{tabular}{|c|c|c|}
\hline
$p$ & $\lnot p$ & $p\vee \lnot p$ \\ \hline
$1$ & $0$ & $1$ \\ \hline
$0$ & $1$ & $1$ \\ \hline
$\frac{1}{2}$ & $\frac{1}{2}$ & $\frac{1}{2}$ \\ \hline
\end{tabular}%
\bigskip
\end{center}

Tye tries to avoid this flaw by saying that a statement is a \emph{%
quasi-tautology} if it has no false substitution instances. So two-valued
tautologies become three-valued quasi-tautologies.

Kleene's three-valued logic is still a precise formalization and having no
tautologies seems a price too high to pay in order to be able to deal in the
above sense with vagueness. Also, according to Keefe, 

\begin{quotation}
[...] the appeal to quasi-tautologies adds nothing: if earning this title is
enough for his [Tye's] purposes, then the fact that $p\vee \lnot p$ also
earns it should be of concern. Moreover, what matters for validity does not
relate to quasitautologies [sic], and assertion depends on sentences being
true not being either true or indefinite, so the role for the notion seems
to be merely one of appeasement. \cite[p. 111]{Keefe}
\end{quotation}

\subsubsection{Fuzzy logics}

Fuzzy logics propose a graded notion of inference. Truth-values range in
degree between $0$ and $1$ in order to capture different degrees of truth.
In this way, the value $0$ is attributed to sentences which are completely
false and the value $1$ to sentences which are completely true. The
remaining sentences are truer than the false sentences, but not as true as
the true ones so they have intermediate logical values according to "how
true" they are. According to Bogenberger

\begin{quote}
In fuzzy logic, the truth of any statement becomes a matter of degree. \cite%
{Bogenberger}
\end{quote}

Fuzzy logic is related to Zadeh's work on fuzzy sets \cite{Zadeh}. A \emph{%
fuzzy set} $A$ on $X$ is characterized by a membership function $f_{A}\left(
x\right) $ with values in the interval $\left[ 0,1\right] $. So, a fuzzy set 
$A$ is a class of objects that allow a continuum of grades of membership.
The membership degree is then the degree to which the sentence '$x$ is a
member of $A$' is true. So, one can interpret the membership degrees of
fuzzy sets as truth degrees of the membership predicate in a suitable
many-valued logic.

Theories of vagueness which recourse to fuzzy logics are advocated most
notably by Machina \cite{Machina} and Smith \cite{Smith}.

According to these theories the notion of heap is a vague one and it may
hold true of given objects only to some (truth) degree. The premises should
be considered partially true to a degree which is quite near to the maximal
degree $1$. This inference has to involve truth degrees for the premises and
has to provide a truth degree for the conclusion in a way that in each step
the truth degree becomes smaller. The sentence '$n$ grains of sand do not
make a heap' tends toward being false for an increasing number of grains.

The problem of saying whether the sentence 'a set of $n$ grains makes a
heap' is true or not is essentially the same as to say that that sentence is
true with a certain (precise) fixed degree. This false precision is perhaps
the main objection to the application of many-valued logics to the sorites
paradox\textit{. }According to Keefe

\begin{quote}
[T]he degree theorist's assignments impose precision in a form that is just
as unacceptable as a classical true/false assignment. In so far as a degree
theory avoids determinacy over whether a is F, the objection here is that it
does so by enforcing determinacy over the degree to which a is F. All
predications of \textquotedblleft is red\textquotedblright \ will receive a
unique, exact value, but it seems inappropriate to associate our vague
predicate \textquotedblleft red\textquotedblright \ with any particular
exact function from objects to degrees of truth. For a start, what could
determine which is the correct function, settling that my coat is red to
degree 0.322 rather than 0.321? \cite[p. 113]{Keefe}
\end{quote}

Also, Urquhart states that

\begin{quote}
One immediate objection which presents itself to [fuzzy logic's] line of
approach is the extremely artificial nature of the attaching of precise
numerical values to sentences like `73 is a large number' or `Picasso's
Guernica is beautiful'. In fact, it seems plausible to say that the nature
of vague predicates precludes attaching precise numerical values just as
much as it precludes attaching precise classical truth values. \cite%
{Urquhart}
\end{quote}

Smith \cite{Smith} tries to solve this problem, suggesting several possible
solutions and concluding that the best answer is to mix fuzzy logic with a
theory called plurivaluationism (not to be confused with supervaluationism%
\footnote{%
Supervaluationism involves only one intended (non-classical) model relevant
to questions concerning meaning and truth, while plurivaluationism allows
that there may be multiple (classical) models.}) called \emph{fuzzy
plurivaluationism}. So, Smith accepts the semantic realism implied by the
Epistemic view, but denies that vague predicates have to refer to a single
bivalent model.

\subsection{Contextualism}

Contextualism\footnote{%
Contextualism is often seen as an argument against philosophical skepticism.
Skepticism claims that we don't actually know what we think we know.
\par
\begin{quote}
But, according to contextualists, the skeptic, in presenting her argument,
manipulates the semantic standards for knowledge, thereby creating a context
in which she can \textit{truthfully} say that we know nothing or very
little. \cite{DeRose}
\end{quote}
} defends that interpretations change over time or according to context.
Such shifts of contexts may occur instantaneously. For instance, at the
beginning of a conversation the context is empty. Then, as the conversation
goes along, these notions are sharpened in such a way that borderline cases
(undecided so far) get assigned to either the extension or the
anti-extension of the vague predicates in question. In fact, borderline
sentences can express something true in one context and something false in
another, so they are context-sensitive. In this way one can disagree about
the truth-values of the propositions expressed by borderline sentences, even
in situations where all the relevant information is available. This view is
most prominently elaborated by Shapiro \cite{Shapiro,Shapiro2008} and
DeRose \cite{DeRose}.

Besides contex-sensitivity Shapiro defines as central the concepts of \emph{%
judgment dependence}, \emph{open texture}, and the \emph{principle of
tolerance}. Judgment dependence means that both the extensions and
anti-extensions for the borderline cases are solely determined by the
decisions of competent speakers. These decisions are put in (and can be
removed from) the \emph{conversational record}. Open texture means that for
a vague predicate $S$ there exists an object $a$ such that a competent
speaker can decide whether $Sa$ holds or not without her competency being
compromised. The principle of tolerance is defined as follows. Suppose that
two objects $a,b$ differ only marginally in the relevant respect on which a
vague predicate $S$ is tolerant. Then if one competently judges $Sa$ to
hold, then $Sb$ also holds.

One reason for skepticism about contextualism is that the problems with
vague expressions seem to remain whether context-sensitivity is taken into
account or not. By taking context into account one can reduce vagueness but
not eliminate completely. Indeed,\ sets with vague boundaries are invariant
to some translations. Take for instance the word 'ugly'. Even if a
particular context is given (and even if one knows a great deal about
another one's ugliness) there is still no reason to suppose that there is a
sharp boundary between what 'ugly' applies to and what it does not.\footnote{The reader interested in the criticism to Shapiro's ideas can consult the review of his book \cite{Shapiro} by Matti Eklund, available online at \url{http://ndpr.nd.edu/news/vagueness-in-context/}. }

Smith \cite{Smith} argues that contextualism should not be seen as a theory
of vagueness in its own right. He claims that this theory is compatible with
all other mentioned theories.

\section{External numbers as a model}
In this section a new approach to the Sorites paradox which takes
advantage of notions and concepts from nonstandard analysis is presented. 
I will start with a brief description of the necessary concepts to reformulate the paradox and model it by means of orders of magnitude. Orders of magnitude are given in the form of the so-called \emph{neutrices} \cite{koudjetivandenberg,dinisberg}.\footnote{The term neutrix was coined by Van der Corput in \cite{VDC} referring to groups of functions, with the intention of creating a mathematical tool that would enable a rigorous \emph{ars negligendi} .}
I want to emphasize that the present response models only a specific type of vagueness
(of the type Sorites) and therefore is not intended as a theory for
vagueness in general. Also, I am not by any means claiming that other theories are
without value. For instance, the fuzzy logic approach has been quite
successful in solving vagueness related to traffic and transportation
processes (see for example \cite{Teodorovic,Bogenberger,Zimmermann} for other examples of applications of fuzzy set theory).
According to Teodorovi\'{c} \cite{Teodorovic}

\begin{quote}
[...] a wide range of traffic and transportation engineering parameters are
characterized by uncertainty, subjectivity, imprecision and ambiguity. Human
operators, dispatchers, drivers and passengers use this subjective knowledge
or linguistic information on a daily basis when making decisions.
\end{quote}

Also,

\begin{quote}
The results obtained show that fuzzy set theory and fuzzy logic present a
promising mathematical approach to model complex traffic and transportation
processes [...]
\end{quote}

However the fuzzy logic approach is also not without fault as model of
imprecision, because it ultimately recourses to precise intervals to model
imprecise situations. Moreover, it does not work with the actual error but
only with an upper bound of the error. Those faults will be corrected with the present proposal.

\subsection{Tools from Nonstandard Analysis} \label{EN}
An ``economical'' version  of Nonstandard Analysis due to Nelson \cite[Chapter 4]{REPT} (see also \cite{dinis}) which is enough for the current purposes is presented below. Add to the language of conventional mathematics a new predicate $\st$. One should read `$x$ is standard' for $\st(x)$. A formula is said to be \emph{internal} if it does not involve the predicate $\st$ (i.e. it is a formula of conventional mathematics) and \emph{external} otherwise. Assume the following:

\begin{enumerate}
\item $\st(0)$;
\item $\forall n \in \mathbb{N} (\st(n) \to \st(n+1))$;
\item $\exists \omega (\neg \st(\omega))$.
\end{enumerate}
Assume also the following axiom scheme.
\begin{axiom}
$$(\Phi(0) \wedge \forallst n (\Phi(n)\to \Phi(n+1)))\to \forallst n \, \Phi(n)$$
\end{axiom}
Were $\Phi$ is an arbitrary formula, internal or external, and $\forallst n \, \Phi(n)$ is an abbreviation of $\forall n (\st(n)\to \Phi(n))$. 
The first two axioms state that the natural numbers from conventional mathematics are all standard. Nevertheless, the third axiom states that there exist nonstandard natural numbers. External induction is a form of induction that allows to conclude that some property is true for all standard natural numbers by assuming that it is valid for 0 and that if it is valid for some standard natural number then it is also valid for its successor.  Of course, the usual form of induction is still valid. However, one should be aware that (internal) induction is a principle from conventional mathematics and is therefore only applicable to internal formulas. Let me illustrate this with a simple example. Let $\Phi(n):\equiv \st(n)$. If we could apply internal induction to $\Phi$ the conclusion would be that $\forall n \, \st(n)$, in contradiction with the third assumption. Induction is applicable to subsets of natural numbers, so we are forced to conclude that $S=\{n: \st(n)\}$ is not a set. These kind of classes are sometimes called \emph{external sets}.

Even this rather weak version of nonstandard analysis is enough to define different orders of magnitude.
\begin{definition}
A real number $x$ is said to be:

\begin{enumerate}
\item \emph{limited}, if there exists $\st \left( n\right) \in \mathbb{N}$ ( i.e., $n\in 
\mathbb{N}
\wedge \st \left( n\right) $) such that $\left \vert x\right \vert
\leq n.$

\item \emph{ilimited}, or \emph{infinitely large} if $x$ is not limited.

\item \emph{infinitesimal}, ou\emph{\ infinitely small} if for any $\st \left( n\right) \in \mathbb{N}_{+}$ one has $\left \vert x\right \vert \leq \frac{1}{n}.$

\item \emph{appreciable }if $x$ is limited but not infinitesimal.

\item Two real numbers whose distance is infinitesimal are said to be \emph{infinitely close}.  
\end{enumerate}
\end{definition}

Sometimes it is not necessary to know precisely the value of a  number to know its order of magnitude. This point of view together intuitions and calculations from nonstandard asymptotics \cite{vdbnaa} lead Koudjeti and Van den Berg \cite{koudjetithese,koudjetivandenberg} to introduce neutrices and external numbers. A \emph{neutrix} is an additive convex subgroup of the reals and an \emph{external number} is the algebraic sum of a real number with a
neutrix. In a nonstandard framework, due to the existence of infinitesimals, there are many neutrices
such as $\oslash$, the external set of all infinitesimals, and
$\pounds $, the external set of all limited numbers (numbers bounded
in absolute value by a standard number).
One can view an external number $\alpha=a+A$ as the sum of a real number $a$ with an ``error'', given by a neutrix $A$. In fact, the rules of calculation for external numbers are a sort of "mellowed" version of the common rules of calculation of real numbers.
Indeed, addition and multiplication in the external numbers are defined (with some abuse of notation) as follows.

\begin{definition}
Let $\alpha=a+A$ and $\beta=b+B$ be two external numbers, the\emph{\ sum} and
\emph{product} of $\alpha$ and $\beta$ are defined as follows
\begin{align*}
\alpha+\beta &  =a+b+max\{A,B\}\\
\alpha \cdot \beta &  =ab+max\{aB,bA,AB\}.
\end{align*}
\end{definition}

The operations are well-defined because neutrices, being convex subgroups of the reals, are ordered by inclusion so the maximum of two neutrices is one of them and because the product of a real number and a neutrix is also a neutrix.
Typically, external numbers are bounded but have neither infimum nor supremum and are stable for
some (but not all!) translations, additions and multiplications. It is not difficult to prove that if $A$ is a neutrix, then for all standard $n$ it holds that $nA=A$. In fact, it holds that $cA=A$, for every appreciable real number $c$. However if we let $\omega$ be nonstandard, then $A \subset \omega A$. Also, if we let $\epsilon$ be infinitesimal we have $\epsilon + \oslash=\oslash$ and even $\oslash + \oslash=\oslash$ but $\epsilon \oslash \subset \oslash$. 
In this way, external numbers generate a calculus of propagation of
errors not unlike the calculus of real numbers, allowing for total order and
even for a sort of generalized Dedekind completeness property \cite{koudjetivandenberg, dinisberg,dinisberg2}.  Thus, external numbers seem suitable as models of orders of magnitude or transitions with imprecise boundaries of the Sorites type, with the advantage of being possible to work directly with imprecisions and errors without recourse to upper bounds. Moreover, the external numbers have a rich algebraic structure. The (external) set of external numbers is a commutative regular semigroup for addition and the (external) set of external numbers which are not reduced to neutrices forms a commutative regular semigroup for multiplication. Although the distributive law does not always hold, necessary and sufficient conditions for it to hold were given. Furthermore the structure has no zero divisors and is sufficiently strong to incorporate combinatorial laws such as the binomial law \cite{dinisberg}. 

A set of individual grains may be
modeled by a standard subset of the external set of limited numbers
(positive part of a neutrix) and the set of grains that form a heap may be
modeled by the external set of the infinitely large numbers.

It should also be possible to capture with external sets some modalities,
like the difference between a "good" approximation, allowing to obtain an
adequately precise numerical result in some context, and a "bad", useless,
one. The stability of orders of magnitude under some repeated additions
justifies to model them by (convex) groups of real numbers.

\subsection[A nonstandard point of view]{A nonstandard point of view on
paradoxical forms\label{section nsa point of view}}

I propose to replace the standard forms presented above by the following
forms which involve reasoning with nonstandard methods.

If one replaces mathematical induction by external induction, the reasoning becomes:%
\begin{equation*}
\left \{ 
\begin{array}{c}
\left( Sa_{1}\wedge \forall ^{st}n\left( Sa_{n}\rightarrow Sa_{n+1}\right)
\right) \rightarrow \forall ^{st}n\,Sa_{n} \\ 
\exists \omega \left( \lnot Sa_{\omega }\right)%
\end{array}%
\right.
\end{equation*}

So, if one admits that:

\begin{enumerate}
\item A single grain of wheat is not a heap.

\item If $n$ is a standard number and if a set of $n$ grains of wheat is not
a heap then a set of $n+1$ grains of wheat is also not a heap.
\end{enumerate}

One concludes that in the presence of a standard number of grains of wheat
one does not have a heap. The heap arises when one has a nonstandard number $%
\omega \simeq +\infty $ of grains of wheat.

The conditional form, using nonstandard analysis, becomes the following.

Let $i$ be a standard natural number. If

\begin{center}
$\left \{ 
\begin{array}{c}
Sa_{1} \\ 
Sa_{1}\rightarrow Sa_{2} \\ 
Sa_{2}\rightarrow Sa_{3} \\ 
... \\ 
Sa_{i}\rightarrow Sa_{i+1}%
\end{array}%
\right. $
\end{center}

Then, by \textit{modus ponens, }the conclusion is $Sa_{i}$, for $i$ an
arbitrarily large but (naive) standard number. In nonstandard analysis this is
modeled by allowing \textit{modus ponens }but only a standard (naive) number
of times. One calls "naive" the natural numbers which can be obtained from
zero by the successive addition of one. This corresponds to Reeb's famous
slogan:

\begin{quote}
Les entiers na\"{\i}fs ne remplissent pas $%
\mathbb{N}
$. \cite{dienerreeb}
\end{quote}

I am in fact claiming that the formalization of the predicate 'is not a
heap' should be an external predicate, where not being a heap means to
possess a standard number of grains.

Soritical arguments share with external numbers the fact of being tolerant
to small changes but not tolerant to large changes in relevant aspects. In
fact, with external numbers, using the different orders of magnitude, it is possible to define rigorously what one
means with terms such as 'small changes' or 'large changes'. The fact that
large changes come as the result of the accumulation of small changes is now
a very natural consequence of the theory.

A simple shift from the classical forms presented in Section \ref{section
induction} and in Section \ref{section conditional form} to the forms using
nonstandard concepts presented in Section \ref{section nsa point of view}\
does not solve the problem. A million grains of wheat should form a heap and
yet that is clearly a standard number of grains. However, both these forms
suggest that the set of individual grains may be modeled by the external set
of limited numbers (positive part of a neutrix) and the set of grains that
form a heap may be modeled by the external set of the infinitely large
numbers. Indeed 'precise' objects possess sharp bounds and can be modeled by
standard sets. 'Vague' objects have no clear bounds and should be for this
matter modeled by nonstandard sets which are given by external properties.

As seen in Section \ref{Epistemic}, epistemicists believe in the existence of sharp bounds for vague concepts, claiming that ignorance is somehow inevitable. The current proposal takes the opposite direction. Indeed, the tolerance of vague terms, such as 'heap', to small changes indicates that such terms do not have a sharp, definite bound. By using neutrices to model such terms it should be possible to avoid the paradox and explain the tolerance to small changes.

According to Keefe \cite{Keefe}, degree theories fail to provide an
acceptable account of vagueness and are forced to make an implausible
commitment to a unique numerical assignment for each sentence. Smith \cite%
{Smith}\ argues that an adequate account of vagueness must involve degrees
of truth and that the main objections to this theory may be overcome. His
fuzzy plurivaluationism theory seems overcomplicated for the present approach to the
Sorites paradox. I believe that the problem with the fuzzy logic approach
is the fact that precise numbers are used to model imprecise predicates. On the contrary, with external numbers 
numbers it is possible to work directly with imprecisions and errors without
recourse to upper bounds, for they have neither infimum nor supremum 
and are tolerant to appreciable (but not infinitely
large) imprecisions. Moreover, since external numbers satisfy strong algebraic laws, similar to the ones of the real numbers, those calculations are quite simple to carry on.

A final remark concerns the strength of nonstandard axioms, which may
introduce undesirable consequences of external modelling. As such, within a nonstandard theory, the proposed solution of the Sorites paradox

\begin{equation}
\left \{ 
\begin{array}{c}
st\left( 0\right) \\ 
\forall n\left( st\left( n\right) \rightarrow st\left( n+1\right) \right) \\ 
\exists \omega \left( \lnot st\left( \omega \right) \right)%
\end{array}%
\right.  \label{SolveSorites}
\end{equation}

implies, by the group property of the standard numbers, invariance by
doubling, i.e. 
\begin{equation*}
\forall n\left( st\left( n\right) \rightarrow st\left( 2n\right) \right) .
\end{equation*}%
One easily imagines a soritical context where this is inappropriate.
However, 
\begin{equation}
\exists n\left( st\left( n\right) \wedge \lnot st\left( 2n\right) \right)
\label{doubling axiom}
\end{equation}%
is consistent with (\ref{SolveSorites}). In such a context \{(\ref%
{SolveSorites}), (\ref{doubling axiom})\} might be an acceptable axiom
system indeed, though of course at the price of losing some calculation properties.

\end{document}